\documentstyle[11pt]{article}

\reversemarginpar
\catcode`\@=11\relax
\newwrite\@unused
\def\typeout#1{{\let\protect\string\immediate\write\@unused{#1}}}
\def\@nnil{\@nil}
\def\@empty{}
\def\@psdonoop#1\@@#2#3{}
\def\@psdo#1:=#2\do#3{\edef\@psdotmp{#2}\ifx\@psdotmp\@empty \else
    \expandafter\@psdoloop#2,\@nil,\@nil\@@#1{#3}\fi}
\def\@psdoloop#1,#2,#3\@@#4#5{\def#4{#1}\ifx #4\@nnil \else
       #5\def#4{#2}\ifx #4\@nnil \else#5\@ipsdoloop #3\@@#4{#5}\fi\fi}
\def\@ipsdoloop#1,#2\@@#3#4{\def#3{#1}\ifx #3\@nnil
       \let\@nextwhile=\@psdonoop \else
      #4\relax\let\@nextwhile=\@ipsdoloop\fi\@nextwhile#2\@@#3{#4}}
\def\@tpsdo#1:=#2\do#3{\xdef\@psdotmp{#2}\ifx\@psdotmp\@empty \else
    \@tpsdoloop#2\@nil\@nil\@@#1{#3}\fi}
\def\@tpsdoloop#1#2\@@#3#4{\def#3{#1}\ifx #3\@nnil
       \let\@nextwhile=\@psdonoop \else
      #4\relax\let\@nextwhile=\@tpsdoloop\fi\@nextwhile#2\@@#3{#4}}
\def\psdraft{
    \def\@psdraft{0}
}
\def\psfull{
    \def\@psdraft{100}
} \psfull
\newif\if@prologfile
\newif\if@postlogfile
\newif\if@bbllx
\newif\if@bblly
\newif\if@bburx
\newif\if@bbury
\newif\if@height
\newif\if@width
\newif\if@rheight
\newif\if@rwidth
\newif\if@clip
\def\@p@@sclip#1{\@cliptrue}
\def\@p@@sfile#1{
           \def\@p@sfile{#1}
}
\def\@p@@sfigure#1{\def\@p@sfile{#1}}
\def\@p@@sbbllx#1{
        \@bbllxtrue
        \dimen100=#1
        \edef\@p@sbbllx{\number\dimen100}
}
\def\@p@@sbblly#1{
        \@bbllytrue
        \dimen100=#1
        \edef\@p@sbblly{\number\dimen100}
}
\def\@p@@sbburx#1{
        \@bburxtrue
        \dimen100=#1
        \edef\@p@sbburx{\number\dimen100}
}
\def\@p@@sbbury#1{
        \@bburytrue
        \dimen100=#1
        \edef\@p@sbbury{\number\dimen100}
}
\def\@p@@sheight#1{
        \@heighttrue
        \dimen100=#1
        \edef\@p@sheight{\number\dimen100}
}
\def\@p@@swidth#1{
        \@widthtrue
        \dimen100=#1
        \edef\@p@swidth{\number\dimen100}
}
\def\@p@@srheight#1{
        \@rheighttrue
        \dimen100=#1
        \edef\@p@srheight{\number\dimen100}
}
\def\@p@@srwidth#1{
        \@rwidthtrue
        \dimen100=#1
        \edef\@p@srwidth{\number\dimen100}
}
\def\@p@@sprolog#1{\@prologfiletrue\def\@prologfileval{#1}}
\def\@p@@spostlog#1{\@postlogfiletrue\def\@postlogfileval{#1}}
\def\@cs@name#1{\csname #1\endcsname}
\def\@setparms#1=#2,{\@cs@name{@p@@s#1}{#2}}
\def\ps@init@parms{
        \@bbllxfalse \@bbllyfalse
        \@bburxfalse \@bburyfalse
        \@heightfalse \@widthfalse
        \@rheightfalse \@rwidthfalse
        \def\@p@sbbllx{}\def\@p@sbblly{}
        \def\@p@sbburx{}\def\@p@sbbury{}
        \def\@p@sheight{}\def\@p@swidth{}
        \def\@p@srheight{}\def\@p@srwidth{}
        \def\@p@sfile{}
        \def\@p@scost{10}
        \def\@sc{}
        \@prologfilefalse
        \@postlogfilefalse
        \@clipfalse
}
\def\parse@ps@parms#1{
        \@psdo\@psfiga:=#1\do
           {\expandafter\@setparms\@psfiga,}}
\newif\ifno@bb
\newif\ifnot@eof
\newread\ps@stream
\def\bb@missing{
    \typeout{psfig: searching \@p@sfile \space  for bounding box}
    \openin\ps@stream=\@p@sfile
    \no@bbtrue
    \not@eoftrue
    \catcode`\%=12
    \loop
        \read\ps@stream to \line@in
        \global\toks200=\expandafter{\line@in}
        \ifeof\ps@stream \not@eoffalse \fi
        \@bbtest{\toks200}
        \if@bbmatch\not@eoffalse\expandafter\bb@cull\the\toks200\fi
    \ifnot@eof \repeat
    \catcode`\%=14
} \catcode`\%=12
\newif\if@bbmatch
\def\@bbtest#1{\expandafter\@a@\the#1
\long\def\@a@#1
\long\def\bb@cull#1 #2 #3 #4 #5 {
    \dimen100=#2 bp\edef\@p@sbbllx{\number\dimen100}
    \dimen100=#3 bp\edef\@p@sbblly{\number\dimen100}
    \dimen100=#4 bp\edef\@p@sbburx{\number\dimen100}
    \dimen100=#5 bp\edef\@p@sbbury{\number\dimen100}
    \no@bbfalse
} \catcode`\%=14
\def\compute@bb{
        \no@bbfalse
        \if@bbllx \else \no@bbtrue \fi
        \if@bblly \else \no@bbtrue \fi
        \if@bburx \else \no@bbtrue \fi
        \if@bbury \else \no@bbtrue \fi
        \ifno@bb \bb@missing \fi
        \ifno@bb \typeout{FATAL ERROR: no bb supplied or found}
            \no-bb-error
        \fi
        \count203=\@p@sbburx
        \count204=\@p@sbbury
        \advance\count203 by -\@p@sbbllx
        \advance\count204 by -\@p@sbblly
        \edef\@bbw{\number\count203}
        \edef\@bbh{\number\count204}
}
\def\in@hundreds#1#2#3{\count240=#2 \count241=#3
             \count100=\count240    
             \divide\count100 by \count241
             \count101=\count100
             \multiply\count101 by \count241
             \advance\count240 by -\count101
             \multiply\count240 by 10
             \count101=\count240    
             \divide\count101 by \count241
             \count102=\count101
             \multiply\count102 by \count241
             \advance\count240 by -\count102
             \multiply\count240 by 10
             \count102=\count240    
             \divide\count102 by \count241
             \count200=#1\count205=0
             \count201=\count200
            \multiply\count201 by \count100
            \advance\count205 by \count201
             \count201=\count200
            \divide\count201 by 10
            \multiply\count201 by \count101
            \advance\count205 by \count201
             \count201=\count200
            \divide\count201 by 100
            \multiply\count201 by \count102
            \advance\count205 by \count201
             \edef\@result{\number\count205}
}
\def\compute@wfromh{
        \in@hundreds{\@p@sheight}{\@bbw}{\@bbh}
        \edef\@p@swidth{\@result}
}
\def\compute@hfromw{
        \in@hundreds{\@p@swidth}{\@bbh}{\@bbw}
        \edef\@p@sheight{\@result}
}
\def\compute@handw{
        \if@height
            \if@width
            \else
                \compute@wfromh
            \fi
        \else
            \if@width
                \compute@hfromw
            \else
                \edef\@p@sheight{\@bbh}
                \edef\@p@swidth{\@bbw}
            \fi
        \fi
}
\def\compute@resv{
        \if@rheight \else \edef\@p@srheight{\@p@sheight} \fi
        \if@rwidth \else \edef\@p@srwidth{\@p@swidth} \fi
}
\def\compute@sizes{
    \compute@bb
    \compute@handw
    \compute@resv
}
\def\psfig#1{\vbox {
    \ps@init@parms
    \parse@ps@parms{#1}
    \compute@sizes
    \ifnum\@p@scost<\@psdraft{
        \typeout{psfig: including \@p@sfile \space }
        \special{ps::[begin]    \@p@swidth \space \@p@sheight \space
                \@p@sbbllx \space \@p@sbblly \space
                \@p@sbburx \space \@p@sbbury \space
                startTexFig \space }
        \if@clip{
            \typeout{(clip)}
            \special{ps:: \@p@sbbllx \space \@p@sbblly \space
                \@p@sbburx \space \@p@sbbury \space
                doclip \space }
        }\fi
        \if@prologfile
            \special{ps: plotfile \@prologfileval \space } \fi
        \special{ps: plotfile \@p@sfile \space }
        \if@postlogfile
            \special{ps: plotfile \@postlogfileval \space } \fi
        \special{ps::[end] endTexFig \space }
        \vbox to \@p@srheight true sp{
            \hbox to \@p@srwidth true sp{
                \hfil
            }
        \vfil
        }
    }\else{
        \vbox to \@p@srheight true sp{
        \vss
            \hbox to \@p@srwidth true sp{
                \hss
                \@p@sfile
                \hss
            }
        \vss
        }
    }\fi
}} \catcode`\@=12\relax

\newenvironment{proof}{\medskip{\bf Proof:}}{\medskip{\hfill $\Box$}}
\def\noproof{\unskip{\hfill $\Box$}}
\newcommand{\lemlab}[1]{\label{lemma:#1}}
\newcommand{\theolab}[1]{\label{theo:#1}}
\newcommand{\codelab}[1]{\label{code:#1}}
\newcommand{\alglab}[1]{\label{alg:#1}}
\newcommand{\eqlab}[1]{\label{eq:#1}}
\newcommand{\corlab}[1]{\label{cor:#1}}
\newcommand{\deflab}[1]{\label{def:#1}}
\newcommand{\tablab}[1]{\label{tab:#1}}
\newcommand{\figlab}[1]{\label{fig:#1}}
\newcommand{\seclab}[1]{\label{section:#1}}
\newcommand{\chaplab}[1]{\label{chapter:#1}}
\newcommand{\openlab}[1]{\label{open:#1}}
\newcommand{\conjlab}[1]{\label{conj:#1}}

\newcommand{\lemref}[1]{\ref{lemma:#1}}
\newcommand{\coderef}[1]{\ref{code:#1}}
\newcommand{\theoref}[1]{\ref{theo:#1}}
\newcommand{\algref}[1]{\ref{alg:#1}}
\newcommand{\corref}[1]{\ref{cor:#1}}
\newcommand{\defref}[1]{\ref{def:#1}}
\newcommand{\tabref}[1]{\ref{tab:#1}}
\newcommand{\figref}[1]{\ref{fig:#1}}
\newcommand{\secref}[1]{\ref{section:#1}}
\newcommand{\chapref}[1]{\ref{chapter:#1}}
\newcommand{\openref}[1]{\ref{open:#1}}
\newcommand{\conjref}[1]{\ref{conj:#1}}

\newcommand{\li}{\item}
\newtheorem{theorem}{Theorem}[section]
\newtheorem{lemma}[theorem]{Lemma}
\newtheorem{prop}[theorem]{Proposition}
\newtheorem{cor}[theorem]{Corollary}
\newtheorem{df}[theorem]{Definition}

{\catcode`\@=11
\gdef\setft#1#2#3{%
\def\@oddfoot{
{\setbox0=\hbox{#1} \setbox1=\hbox{#3} \ifdim\wd0>\wd1
\dimen0=\wd0 \box0\hfil#2\hfil\hbox to\dimen0{\hfil\hfil\box1}
\else \dimen0=\wd1 \hbox to\dimen0{\box0\hfil }\hfil#2\hfil\box1
\fi }}} }

\def\complaint#1{}
\def\withcomplaints{
\newcounter{mycomplaints}
\def\complaint##1{\refstepcounter{mycomplaints}%
\ifhmode%
\unskip%
{\dimen1=\baselineskip \divide\dimen1 by 2 %
\raise\dimen1\llap{\tiny -\themycomplaints-}}\fi%
\marginpar{\tiny [\themycomplaints]: ##1}}%
}

\begin{document}
\bibliographystyle{plain}

\title{On the Number of Embeddings of Minimally Rigid Graphs}

\author{\setcounter{footnote}{0}%
\def\thefootnote{\arabic{footnote}}
Ciprian Borcea \footnotemark[1]\hspace{3mm} and Ileana~Streinu
\footnotemark[2]\hspace{3mm}  }

\footnotetext[1]{Department\ of\ Mathematics, Rider University,
Lawrenceville, NJ 08648, USA, borcea@rider.edu}

\footnotetext[2]{Department\ of\ Computer\ Science, Smith College,
Northampton, MA 01063, USA, streinu@cs.smith.edu. Supported by NSF
grants CCR-0105507 and CCR-0138374.}

\date{} 
\maketitle

\begin{abstract}

Rigid frameworks in some Euclidian space are embedded graphs
having a unique local realization (up to Euclidian motions) for
the given edge lengths, although globally they may have several.
We study  the number of distinct planar embeddings of minimally
rigid graphs with $n$ vertices. We show that, modulo planar rigid
motions, this number is at most ${{2n-4}\choose {n-2}} \approx
4^n$. We also exhibit several families which realize lower bounds
of the order of $2^n$, $2.21^n$ and $2.88^n$.

For the upper bound we use techniques from complex algebraic
geometry, based on the (projective) Cayley-Menger variety
$CM^{2,n}(C)\subset P_{{{n}\choose {2}}-1}(C)$ over the complex
numbers $C$. In this context, point configurations are represented
by coordinates given by squared distances between all pairs of
points. Sectioning the variety with $2n-4$ hyperplanes yields at
most $deg(CM^{2,n})$ zero-dimensional components, and one finds
this degree to be   $D^{2,n}=\frac{1}{2}{{2n-4}\choose {n-2}}$.
The lower bounds are related to inductive constructions of
minimally rigid graphs via Henneberg sequences.

The same approach works in higher dimensions. In particular we
show that it leads to an upper bound of $2 D^{3,n}=
{\frac{2^{n-3}}{n-2}}{{n-6}\choose{n-3}}$ for the number of
spatial embeddings with generic edge lengths of the $1$-skeleton
of a simplicial polyhedron, up to rigid motions.

\end{abstract}

\section{Introduction}
\label{introduction}

In this paper we are concerned with {\em graph embeddings subject
to edge lengths constraints}. We use {\it embedding} (or {\it
realization}) in the extended sense, which allows some vertices to
coincide and some edges to cross. For a given graph and for a
given set of edge lengths, a natural question to ask is: {\it how
many embeddings in $R^d$ are there?}

Obviously, for a fixed dimension $d$, some graphs  have a
continuum of embeddings, or may have no embedding at all for
particular choices of edge lengths.

We consider {\it 
minimally rigid graphs on $n$ vertices in dimension $d$} (which
have $dn-{{d+1}\choose 2}$ edges and unique local realizations for
generic choices of edge lengths), with particular regard to
dimensions $2$ and $3$.

Our results give a general upper bound in arbitrary dimension,
which is of the order of $2^{dn}$ for fixed $d$ and $n$
sufficiently large. We also exhibit a family of graphs inducing a
lower bound of the order of $2.88^n$ in dimension $2$.

\medskip
\noindent {\bf Historical Perspective.} {\it Distance geometry}
relies on foundational work of Cayley on configurations of $n$
points in Euclidean d-space. As elaborated by Menger (see
\cite{berger}, p. 237), this led to conditions characterizing
systems of ${n \choose 2}$ positive reals that arise as pairwise
squared distances between $n$ points in $R^d$. See also
\cite{blumenthal} and \cite{deza}. In dimension $3$, distance
geometry has been used for the study of molecular conformation in
chemistry (\cite{crippen}, \cite{crippenhavel},
\cite{hendrickson2}). Indeed, inter-atomic distance information
can be obtained from nuclear magnetic resonance spectra of a
molecule. Solving the graph embedding problem determines
coordinates for the atoms, and hence the $3$-dimensional shape of
the molecule. Other applications include surveying and satellite
ranging.

Our investigation is also related to classical studies in the
Kinematics of mechanical linkages, in particular the problem of
tracing algebraic curves. Wunderlich (\cite{wunderlich}) gives an
interesting family of generalized planar coupler curves with
degree growing exponentially in the number of links. Related work
was also  done in combinatorial Rigidity Theory (\cite{connelly},
\cite{walterS}, \cite{servatius}). Rigid frameworks are embedded
graphs having a unique local realization for the given edge
lengths. But globally there may be several realizations. Using
special combinations of graphs and edge lengths, Saxe \cite{saxe}
has shown that it is NP-hard to solve the graph embedding problem
in dimension two, as well as to determine whether it has a unique
solution. Under the assumption of genericity, Hendrickson
\cite{hendrickson} studied conditions on rigid frameworks that
guarantee a unique global realization.

\medskip
\noindent {\bf Minimally Rigid Graphs.} Hendrickson's graphs, as
well as the complete graphs implicit in the Cayley-Menger
conditions, are highly redundant in terms of distance
dependencies. This makes the problem of reconstructing point sets
from distances very challenging, as small perturbation errors in
the input data will render the problem infeasible. In this paper
we go to the lower end of the dependency spectrum and study the
number of
embeddings of {\it 
minimally rigid graphs in dimension $d$}.

In dimension two, these graphs possess a simple combinatorial
characterization \cite{laman}, and, for {\em generic} edge
lengths, will have a {\em discrete} realization space. We'll refer
to these graphs as {\it Laman graphs}. In dimension three, one
does not have a comparable result. However, as observed by Gluck
\cite{gluck}, a combination of arguments going back to Cauchy,
Steinitz, Dehn, Weyl and Alexandrov ensures that $1$-skeleta of
simplicial polyhedra
are 
minimally rigid graphs in dimension three.

\medskip
\noindent {\bf Techniques and Results.} Our approach emphasises
the possibility of exploiting the {\em algebraic} character of the
prevalent relations in distance geometry by rephrasing the
problems in terms of associated {\em complex projective
varieties}.

Indeed, our upper bounds are derived from {\em degree}
computations for the complex projective varieties $CM^{d,n}$
induced by the Cayley-Menger determinant equalities for squared
distances of $n$ points in $R^d$. In \cite{borcea} these {\em
Cayley-Menger varieties} are identified with classical {\em
determinantal varieties} of symmetric forms.

\medskip
For planar realizations of Laman graphs we obtain:
\begin{theorem}
\label{twodim} Given a generic choice of edge lengths, a Laman
graph with $n$ vertices has at most ${{2n-4}\choose {n-2}}$ planar
embeddings, up to rigid motions.
\end{theorem}

\medskip
In the direction of lower bounds, we exhibit several families
which realize bounds of the order of $2^n$, $2.21^n$ and $2.88^n$.

For spatial realizations of $1$-skeleta of simplicial polyhedra we
obtain:

\begin{theorem}
\label{threedim} Given a generic choice of edge lengths, the
$1$-skeleton of a simplicial complex which has $n$ vertices and is
topologically a $2$-dimensional sphere \ \ has \  at \  most
${\frac{2^{n-3}}{n-2}}{{n-6}\choose{n-3}}$ spatial embeddings, up
to rigid motions.
\end{theorem}

\medskip
In general, we have:

\begin{theorem}
\label{ddim} Let $d$ be a given dimension. Let $\cal G$  be a
class of connected graphs with $n\geq d+1$ vertices and
$m=dn-{{d+1}\choose 2}$ edges such that each graph allows an
infinitesimally rigid realization in $R^d$.

Then, for a generic choice of edge lengths, each graph in $\cal G$
has a finite number of embeddings in $R^d$ and this number is
bounded from above by $2 D^{d,n}$. Here $D^{d,n}$ stands for the
degree of the Cayley-Menger variety $CM^{d,n}(C) \subset
P_{{n\choose 2}-1}(C)$ and is given by the formula:

$$ D^{d,n} = \prod_{k=0}^{n-d-2} \frac {{{n-1+k} \choose {n-d-1-k}}}{{{2k+1}\choose {k}}} $$

For $d$ constant and $n$ sufficiently large, $D^{d,n} \approx
2^{dn}$.

\end{theorem}

\medskip
It is worth noticing  that known (but general) techniques based on
{\it real algebraic geometry} (Oleinik-Petrovskii-Milnor-Thom
bounds) yield less good upper bounds.

\medskip
\noindent {\bf Future Perspectives and Open Questions.} Trying to
bridge the still significant gap between the upper and lower
bounds can be approached from both directions. The upper bound is
presumably  high, as it counts not only ``realistic'' solutions,
but also real solutions which do not represent configurations  and
complex solutions in the Cayley-Menger variety.

On the other hand, in showing that Laman graphs which require
Henneberg II constructions may have more than $c^n, c>2$ planar
embeddings, we give evidence that the lower bound may possibly be
raised.

Obviously, an important problem for $d\geq 3$ is to identify and
characterize classes of minimally rigid graphs.

\section{Definitions and Preliminaries}
\label{definitions}

In this section we define the necessary concepts: configuration
space, rigid, infinitesimally rigid, $d$-minimally rigid graph. We
also formulate the problem and, for the sake of comparison,
 use standard results from Real Algebraic Geometry to derive  upper bounds that
we'll later improve.

\medskip
\noindent {\bf Rigidity Theoretic Definitions.} \label{rigidity}
Let $G$ be a graph $G=(V,E)$, $V=\{1, \cdots, n\}$, $m=|E|$. A
{\it framework} $(G,L)$ is a graph $G$ together with a set
$L=\{l_{ij} | ij\in E\}$ of positive numbers $l_{ij}> 0$
interpreted as {\it edge lengths} associated to the edges. A {\it
realization} or {\it embedding} $G(P)$ of $(G,L)$ in some space
$R^d$ is given by a mapping of the vertices onto a set of points
$P=\{p_1, \cdots, p_n\}\subset R^d$ such that $l_{ij}$ equals the
Euclidian distance between the two points $p_i$ and $p_j$. Note
that edges may cross and vertices may coincide.

\medskip

{\bf Notation.} The squared distance (resp. the distance) between
two points labeled $i$ and $j$ will be denoted by $d_{ij}$ (resp.
$l_{ij}$) when given a priori, and by $x_{ij}$ (squared distance
only) when unknown.

\medskip

The {\it realization} or {\it configuration space of $(G,L)$ in
$R^d$} is the set of all possible realizations,
modulo the ${{d+1}\choose 2}$ dimensional group of rigid motions
(translations and orthogonal transformations).

Any graph has some embedding, for some edge lengths, but for given
values of $L$ the configuration space may be empty, finite or
higher dimensional. A realization $G(P)$ is {\bf rigid} when it
cannot be deformed continuously into another (non-congruent)
realization of the same framework. Otherwise, the realization is
just one point of a higher dimensional component of the
realization space, and the embedded graph is called in this case a
{\it mechanism}, as it is {\it flexible}. The dimension of the
local component to which the realization belongs is its {\it
number of degrees of freedom}. A {\it one-degree-of-freedom
mechanism} (abbreviated as {\it 1DOF mechanism}) is an embedded
framework whose component is
one-dimensional (i.e. a curve). 

The {\bf rigidity matrix} $A$ associated to a realization $G(P)$
is the following $m\times dn$ matrix. The rows are indexed by the
$m$ edges $ij\in E$ and the $dn$ columns are grouped into
$d$-tuples corresponding to the $d$ coordinates of a point indexed
by $i$, $i=1, \cdots, n$. The $ij$th row has $0$ entries
everywhere, except in the $i$th and $j$th $d$-tuples of columns,
where the entries are $p_i-p_j$, resp. $p_j-p_i$.

$$
\begin{array}{ll}
\ & \begin{array}{lllllll} \ \ 1 & \ \ \cdots & \ \ \ \  i &\  \ \
\ \cdots & \ \ \ \ j & \ \ \ \cdots & \ \ n
\end{array} \\

ij & \left(
\begin{array}{lllllll}
\ & \cdots & \ & \ & \ & \ & \ \\
0 & \cdots & p_i -p_j & \cdots & p_j -p_i & \cdots & \ 0 \\
\cdots & \ & \ & \ & \ & \
\end{array}
\right)
\end{array}
$$

A solution $v\in R^{dn}$ of the linear system $Av=0$ is called an
{\it infinitesimal motion}. It is {\it trivial} if arising from an
infinitesimal rigid transformation of $R^d$. The space of trivial
infinitesimal motions is a ${{d+1}\choose 2}$-dimensional linear
subspace of the space of solutions of $Av=0$. Hence the rank of
the rigidity matrix is at most $dn-{{d+1}\choose 2}$. It is known
that it depends on both the combinatorial structure of the graph
$G$ and on the points of the embedding $P$.

A realization is {\bf infinitesimally rigid} if the rank of its
associated rigidity matrix is precisely $dn-{{d+1}\choose 2}$.
This means that the only solutions $v\in R^{dn}$ of the system
$Av=0$ are the trivial infinitesimal motions. An infinitesimally
rigid configuration is rigid, but the opposite is not true: there
exist rigid embeddings which are not infinitesimally rigid.

\begin{df}
For a given dimension $d$, we say that a graph with $n$ vertices
is {\bf $d$-minimally rigid} or {\bf minimally rigid in dimension
$d$} if it has $m=nd - {{d+1}\choose 2}$ edges and {\it allows an
infinitesimally rigid realization in $R^d$}. These graphs are also
called {\bf $d$-isostatic}, see \cite{servatius}.
\end{df}

{\bf Remark.} As a maximal rank condition, infinitesimal rigidity
will hold true on the complement of
an algebraic subvariety if holding true at one point. Consequently, a 
$d$-minimally rigid graph, as defined above, will have only
infinitesimally rigid realizations in $R^d$, provided
$d_{ij}=l^2_{ij}, ij\in E$ are chosen in the complement of a
certain algebraic subvariety (i.e. {\it generically} as the term
is used in algebraic geometry).

\medskip
\noindent {\bf DIMENSION TWO.} Laman \cite{laman} gave a complete
characterization
of 
$2$-minimally rigid graphs: we will call them {\it Laman graphs}.
They have $2n-3$ edges, and each subset of $k$ vertices spans at
most $2k-3$ edges.

\medskip
\noindent {\bf Henneberg Sequences.} \label{henneberg} Laman
graphs have many elegant combinatorial properties. The so-called
{\it Henneberg sequences} (see \cite{henneberg}, \cite{walter})
are used to construct them in an inductive fashion.  A Henneberg
sequence for a graph $G$ is a sequence $G_3, G_4, \cdots, G_n$ of
Laman graphs on $3, 4, \cdots, n$ vertices, such that: $G_3$ is a
triangle, $G_n=G$ and each graph $G_{i+1}$ is obtained from the
previous one $G_i$ via one of two types of {\it steps}: type I and
type II. A Henneberg step of type I adds a new vertex and two new
edges connecting this vertex to two arbitrary vertices of $G_i$. A
Henneberg step of type II adds a new vertex and three new edges,
and removes an old edge, more precisely: the three new edges must
connect the new vertex to three old vertices such that at least
two of them are joined via an edge; that edge will be removed. See
Figure \figref{henne}.

\begin{figure}[htbp]
\vspace{-0.2in}
\begin{center}
\ \psfig{file=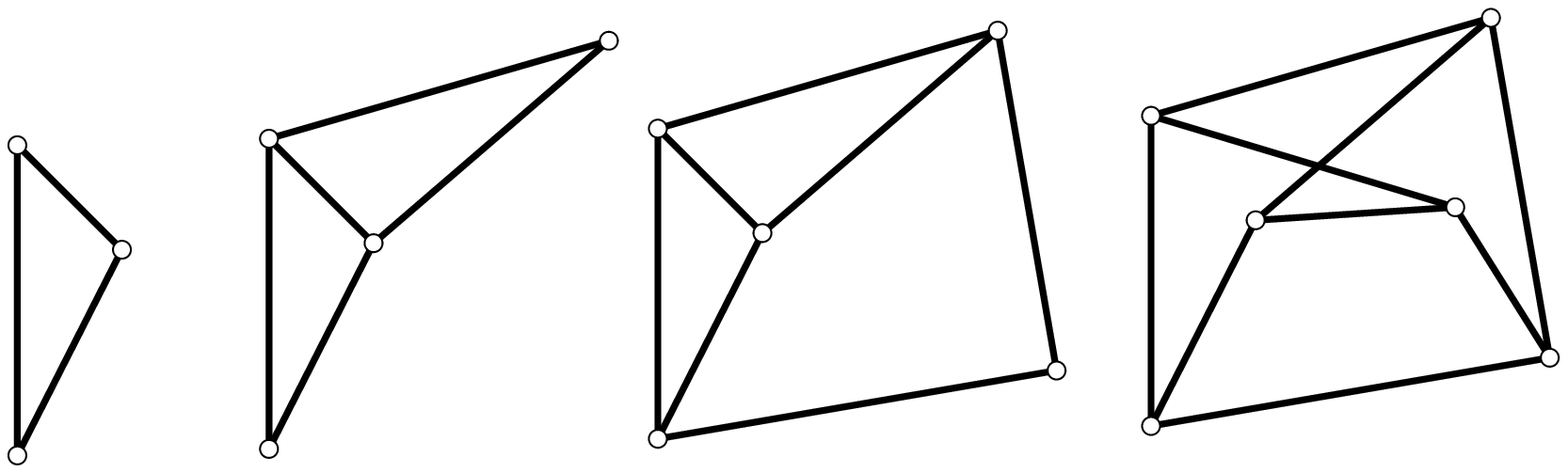,width=3.2in}
\psfig{file=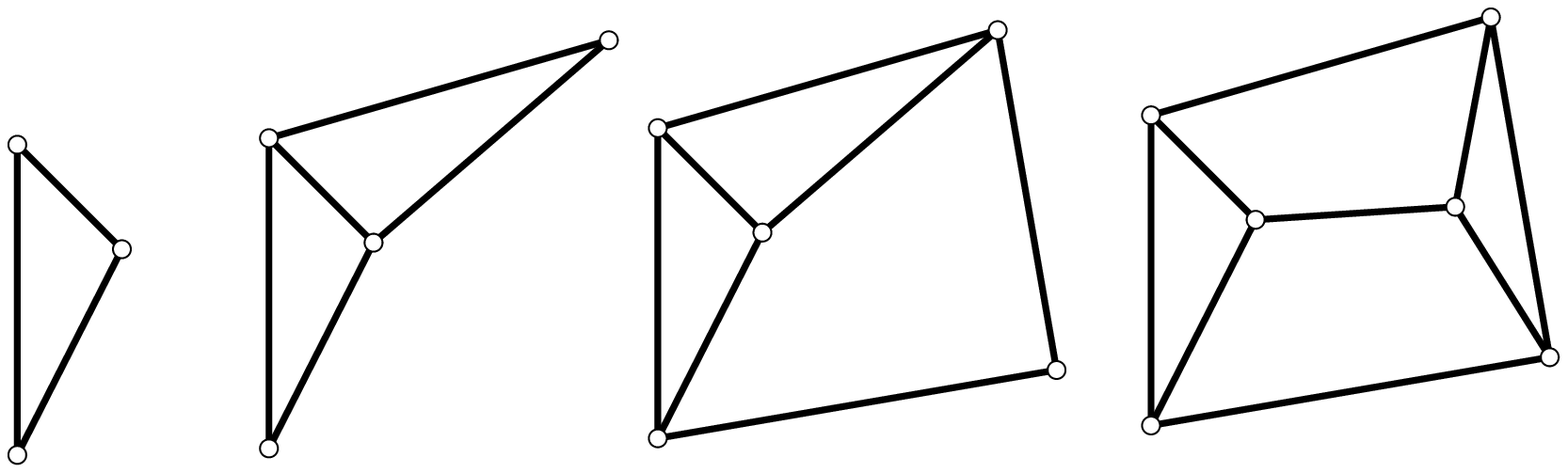,width=3.2in}
\end{center}
\vspace{-0.2in} \caption{
Henneberg constructions for two graphs on $6$ vertices, using type
I and type II steps. Left: $K_{3,3}$. Right: the Desargues graph.}
\figlab{henne}
\end{figure}


Not all Laman graphs can be constructed using only type I steps.
We call {\it Henneberg I graphs} those which can, and {\it
Henneberg II} those which cannot. In particular, triangulations of
planar polygons are Henneberg I graphs. The smallest examples of
Henneberg II graphs are shown in Figure \figref{henne}. The one on
the right, the so-called {\it Desargues framework}, borrows its
name from the classical Desargues configuration (two triangles in
perspective) from plane projective geometry, which it induces for
some special edge lengths.  As an aside, we note that these
special embeddings, illustrated in the last two cases of Fig.
\figref{desnr}, are not infinitesimally rigid, or may even be
flexible. Such embeddings are not {\it generic}, and therefore not
treated in the general discussion of this paper. The first case in
Fig. \figref{desnr} is an example of a {\it generic} embedding.

\begin{figure}[htbp]
\vspace{-0.2in}
\begin{center}
\ \psfig{file=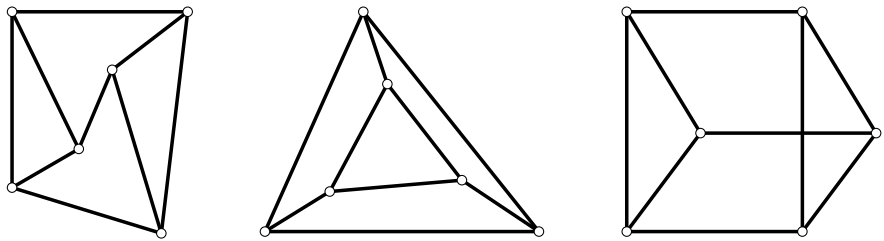,height=1.4in}
\end{center}
\vspace{-0.2in}
\caption{
The Desargues framework and three types of its possible
embeddings: infinitesimally rigid, rigid but not infinitesimally
rigid and flexible.} \figlab{desnr}
\end{figure}

\medskip
\noindent {\bf DIMENSION THREE.} \label{dim3}
A combinatorial characterization of 
minimally rigid graphs in dimension three or higher is not
yet known. However, there is an important class of graphs which are 
minimally rigid in dimension $3$, namely the $1$-skeleta of
simplicial convex polyhedra. By a theorem of Steinitz
\cite{steinitz}, these are the $3$-connected maximally planar
graphs (planar triangulations). The realization of such a graph as
a $1$-skeleton of a simplicial convex polyhedron is
infinitesimally rigid, see Gluck \cite{gluck}.

\medskip
\noindent {\bf HIGHER DIMENSIONS.} \label{dimd} Here is a simple
argument showing that there exist generic minimally rigid graphs
in any dimension $d$. Consider the class of graphs constructed
inductively as follows. Start with the $1$-skeleton of a
$d$-simplex (the complete graph on $d$ vertices). At each step
afterwards, add a new vertex and $d$ edges connecting it to $d$ of
the previously constructed vertices. It is straightforward to
construct an infinitesimally rigid embedding of this graph, by
embedding the vertices on a set of points in general position (an
induction on the construction shows that the rank of the rigidity
matrix is as expected).

It would be interesting to find other classes of minimally rigid
graphs with simple descriptions.

\medskip
\noindent {\bf The Embedding Problem.} \label{embedding} Given a
framework $(G,L)$, one is interested in understanding its
configuration space; in particular, topological invariants such as
the number of connected components and their dimensions.


The configuration space is the orbit space of the group of rigid
motions acting on the affine algebraic variety defined in $R^{dn}$
by the equations:
\begin{equation}
\label{eq1} |p_i-p_j|^2=d_{ij}, \ ij\in E
\end{equation}

Since the  group is connected, the number of connected components
can be counted on the algebraic variety itself.

\medskip
\noindent {\bf Bounds from Real Semi-Algebraic Geometry.}
\label{semialgebraic} Consider an algebraic system of $s$
equations of degree $d$ in $k$ variables $P_l[x_1, \cdots, x_k]=0,
l=1, \cdots, s$, degree $P_l \leq d$. Its real solutions form a
{\it real algebraic set}.

The following result from real algebraic 
geometry bounds the number of connected components or, more
generally, the Betti numbers of an algebraic set.

\medskip
{\bf Oleinik-Petrovskii-Milnor-Thom Theorem} (\cite{oleinik},
\cite{milnor}, \cite{thom}, see also \cite{bpr2}): for equations
of degree at most $d$, in $k$ variables, the sum of the Betti
numbers (and hence the number of connected components) is at most
$d(2d-1)^{k-1}$.

\medskip
Applying the theorem to a minimally rigid graph in dimension $d$:
the system has $m=dn-{{d+1}\choose 2}$ equations of degree $2$ in
$dn$ variables. The Oleinik-Petrovsky-Milnor-Thom bound implies at
most $2\times 3^{dn-1} \approx (3^d)^n$ components (i.e.
realizations). In particular, in dimension $2$ the bound is of the
order  $9^n$.

\medskip
We will obtain better upper bounds on the number of realizations
by using a different algebraic object, the Cayley-Menger variety.

\section{Cayley-Menger Varieties}
\label{cmvariety}

A set of ${n \choose 2}$ positive numbers $d_{ij}=l_{ij}^2$,
$i,j=1, \cdots, n, i<j$, has to satisfy certain algebraic
conditions to be the set of squared Euclidian distances between
$n$ points in $R^d$. Some of these conditions are equalities (due
to Cayley) and express the fact that the rank of the following
matrix is at most $d+2$:

$$
\left(
\begin{array}{lllll}
0 & 1 & 1 & \cdots & 1\\
1 & 0 & d_{12} & \cdots & d_{1n} \\
1 & \  & \cdots & \ & \ \\
1 & d_{n1} & d_{n2} & \cdots & 0
\end{array}
\right)
$$

The other conditions, due to Menger, are inequalities.

If we restrict our attention only to the equalities, we are led to
consider a {\em complex projective variety} $CM^{d,n}(C)\subset
P_{{{n}\choose {2}}-1}(C)$, which is called a {\it Cayley-Menger
variety} in \cite{borcea}. Here $C$ stands for the field of
complex numbers, to emphasize that we refer to the complex and not
the real variety. It is defined by setting to zero all the
$(d+3)\times (d+3)$ minors of the following {\em symmetric}
matrix:

$$
\left(
\begin{array}{lllll}
0 & 1 & 1 & \cdots & 1\\
1 & 0 & x_{12} & \cdots & x_{1n} \\
1 & \  & \cdots & \ & \ \\
1 & x_{n1} & x_{n2} & \cdots & 0
\end{array}
\right)
$$

\noindent where $(x_{ij})_{1\leq i < j\leq n}$ are {\em complex
homogeneous coordinates} in the projective space $ P_{{{n}\choose
{2}}-1}(C)$.

We arrive in this manner to a {complex-projective} formulation of
our Embedding Problem: {\it Solve the system consisting of the
Cayley determinantal equations, plus the linear conditions setting
the graph edges to prescribed lengths.} The latter conditions have
the form:

$$   \frac{x_{ij}}{x_{kl}}=\frac{d_{ij}}{d_{kl}}, \ \ \mbox{that is:} \  \
d_{kl}x_{ij}=d_{ij}x_{kl}, \ \ \mbox{for} \ ij,kl\in E $$

\noindent which (for $m=|E|=dn-{{d+1}\choose 2}$ edges) amounts to
$dn-{{d+1}\choose 2} -1$ independent {\em hyperplane sections}.

\medskip
As we are about to see, this is precisely the (complex) dimension
of the Cayley-Menger variety $CM^{d,n}(C)$. Hence, by a general
result in complex projective geometry, we have:

\begin{prop}
A codimension $dn-{{d+1}\choose 2} -1$ linear section of the
Cayley-Menger variety $CM^{d,n}(C)$ has  at   most
$deg(CM^{d,n}(C))$ isolated points.
\end{prop}

\noindent Here, $deg(CM^{d,n}(C))$ stands for the {\em degree} of
our variety, which is defined as the number of (isolated) points
in a {\em generic} codimension $dn-{{d+1}\choose 2} -1$ linear
section.

\medskip \noindent
{\bf Remark:} Even with generic choices of edge lengths
$(l_{ij})_{ij\in E}$,
the linear sections corresponding to 
$d$-minimally rigid  graphs remain very peculiar, and are never
generic in the sense required in the definition of the degree.
However, we know that a generic choice of edge lengths ensures
that all configuration solutions, (which we also call
``realistic'', meaning both real {\em and} satisfying Menger's
inequalities) are isolated points, because the resulting framework
is infinitesimally rigid.

\begin{cor}
\label{cor1}
For a generic choice of edge lengths, a 
$d$-minimally rigid \ \ graph \  \  with $n$ vertices has \ \ at \
\ most $2deg(CM^{d,n}(C))$ distinct planar embeddings, up to rigid
motions.
\end{cor}

Indeed, there will be  twice as many solutions to the original
embedding problem than ``realistic'' solutions in the sectioning
of the Cayley-Menger variety, because we count as distinct two
realizations which are one the reflection of the other (modulo
rigid transformation in the plane), while the Cayley-Menger
approach automatically identifies congruent configurations. (The
projective equivalence, which identifies similar configurations is
no issue, since one can always consider a `first' edge as the unit
of length.)

\section{The Upper Bounds}
\label{upper}

We begin with a ``naive'' count for dimension: $n$ points in $R^d$
require $dn$ parameters; equivalence under Euclidean motions and
rescaling cuts down ${{d+1}\choose 2}+1$ parameters, so that the
{\em  configuration space} for $n$ (ordered) points in the plane,
modulo congruence and similarity, should be $(dn-{{d+1}\choose
2}-1)$-dimensional. The complex version $CM^{d,n}(C)$ (which is,
technically, the Zariski-closure \ of \ the \  configuration space
in $P_{{{n}\choose {2}}-1}(C)$), \ should have as many complex
dimensions.

\medskip
In order to be precise, we first relate {\em Cayley coordinates}
and {\em Gram coordinates} for \  a \  configuration \  of \  $n$
\  points $p_1,...,p_n\in R^d$.

\medskip
The {\em Cayley coordinates} are:
$$ x_{ij}=x_{ij}(p)=|p_i-p_j|^2=\langle p_i-p_j,p_i-p_j\rangle, \ \ \ 1\leq i,j \leq n $$
\noindent with $\langle\ ,\ \rangle$ denoting the usual inner
product in $R^d$.

\medskip
For {\em Gram coordinates}, we have to {\bf choose} one point as
origin, and {\bf we make the choice} $p_1=0$. Then we set:
$$ y_{ij}=y_{ij}(p)=<p_i-p_1,p_j-p_1>=<p_i,p_j>, \ \ \ 2\leq i,j \leq n $$
The relation is simply the {\em cosine theorem}:
$$ y_{ij}=\frac{1}{2}(x_{1i}+x_{1j}-x_{ij}), \ \ \ 2\leq i, j\leq n $$
Normally, one looks at the Gram coordinates as arranged in a {\em
symmetric $(n-1)\times (n-1)$ matrix} $Y$ with entries
$y_{ij}=y_{ji}$, while we've seen above the Cayley coordinates
arranged in a {\em bordered symmetric $(n+1)\times (n+1)$ matrix},
which we denote by $X$.

\medskip
It is a simple exercise (cf. \cite{berger} or \cite{borcea}) to
observe the rank relation:
$$ rk(X)=2+rk(Y) $$
Since for configurations in $R^d$ the Gram matrix has obviously
rank at most $d$, we see the reason for defining the Cayler-Menger
variety by the vanishing of all $(d+3)\times (d+3)$ minors in $X$.
In fact, we have:

\begin{prop}
The passage from Cayley coordinates to Gram coordinates identifies
the Cayley-Menger variety $CM^{d,n}(C)$ with the determinantal
variety defined (projectively) by all complex (non-trivial)
symmetric $(n-1)\times (n-1)$ matrices of rank at most $d$.
\end{prop}

This brings the matter into the territory of classical algebraic
geometry, and we may simply refer to \cite{Har} and \cite{Ful}.
The degree computation  goes back to Giambelli
 (cf. \cite{Gia} \cite{H-T} \cite{JLP}).

The result converts our Corollary \ref{cor1} into:

\begin{theorem}
Let $d$ be given and assume a generic choice of edge lengths.
A 
$d$-minimally rigid graph with $n$ vertices has at most
$$2 D^{d,n} = 2 \prod_{k=0}^{n-d-2} \frac {{{n-1+k} \choose {n-d-1-k}}}{{{2k+1}\choose {k}}} $$
embeddings in $R^d$, up to rigid motions.
\end{theorem}

This is our Theorem \ref{ddim} in the Introduction. A direct
calculation shows that this number is $O(2^{dn})$ for $d$ fixed
and $n$ sufficiently large.

\medskip

For dimension two we have:

\begin{theorem}
Given a generic choice of edge lengths, a Laman graph with $n$
vertices has at most  $2deg(CM^{2,n}(C))={{2n-4}\choose {n-2}}$
planar embeddings, up to rigid motions.
\end{theorem}

A direct computation of this particular degree is given in
\cite{borcea}.

In dimension three, we obtain:

\begin{theorem}
Given a generic choice of edge lengths, the $1$-skeleton of a
simplicial convex polyhedron with $n$ vertices has at most
$2deg(CM^{3,n}(C))=\frac{2^{n-3}}{n-2}{{3n-6}\choose {n-3}}$
embeddings in $R^3$, up to rigid motions.
\end{theorem}

\begin{figure}[htbp]
\vspace{-0.3in}
\begin{center}
\ \psfig{file=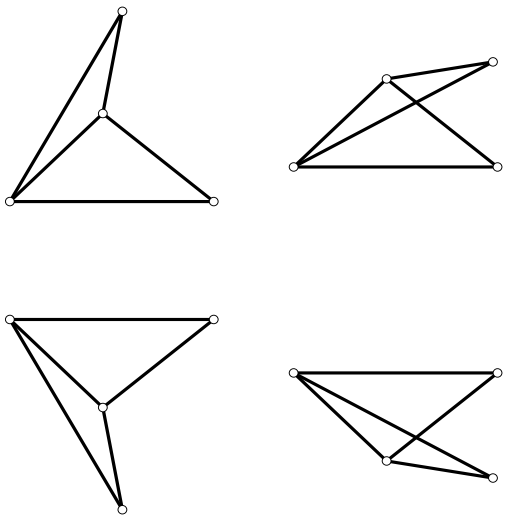,height=2.8in,width=2.8in}
\end{center}
\vspace{-0.3in} \caption{
Triangulations have $2^{n-2}$ embeddings. An example for $n=4$.}
\figlab{triang} \vspace{-0.2in}
\end{figure}

\section{Lower Bounds}
\label{lower}

The simplest family of examples yielding an exponential lower
bound is given by triangulations of simple polygons. Indeed, any
triangle may be flipped over an adjacent internal diagonal, giving
a total of exactly $2^{n-3}$ embeddings. Taking into account the
mirror image of the whole graph doubles the number of
possibilities to $2^{n-2}$. See Fig. \figref{triang}.

Our goal in this section is to improve this trivial bound. We do
so by refining the analysis of minimally rigid graphs via
Henneberg constructions. We show first that Henneberg I graphs
have at most as many embeddings as triangulations have, and for
particular choices of edge lengths this number is achieved. Then
we improve the lower bound twice, using two types of special
iterations of Desargues graphs.

\medskip
\noindent {\bf Bounds on the number of embeddings for Henneberg I
graphs.} Let's first analyze the Henneberg constructions for
triangulations of polygons. Starting from any triangle, at each
step we add a new vertex and two new edges connecting it to two
old vertices {\it which are already joined by an edge}, hence at a
fixed distance. When the edge lengths are given, this construction
is easily carried out using ruler and compass. The newly added
vertex can be chosen to be one of the two intersections of two
circles with radii the edge lengths of the newly added edges
(which are known to intersect, as the edge lengths are computed
from an existing triangulation.) Any choice at step $i$ is valid
and does not affect the number of further options.

\begin{figure}[htbp]
\vspace{-0.3in}
\begin{center}
\ \psfig{file=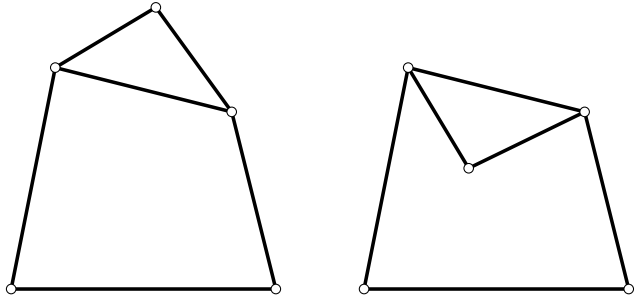,height=1.8in}
\end{center}
\vspace{-0.3in}
\caption{
The $4$-bar mechanism obtained from the Desargues graph, in two
embeddings of the coupler triangle.} \figlab{des1}
\end{figure}

We now see that this analysis depends only on the fact that, at
step $i$ in the Henneberg I construction, the two vertices to be
joined by the newly added edges to the new vertex are themselves
connected by an edge, hence at fixed distance. This amounts to
``adding a triangle'' to the previously constructed graph.
Therefore, for a realizable choice of edge lengths (which always
exists) we get:

\begin{lemma}
Any realizable Henneberg I graph obtained by ``adding triangles''
has exactly $2^{n-2}$ embeddings.
\end{lemma}

In the general Henneberg I step, the two added edges may not form
a triangle with an old edge. Therefore the distances between
relevant points (such as those used to add two new diagonals at
the next step) may vary between embeddings constructed up to step
$i$, and this may render the next step infeasible (i.e. having two
complex, instead of two real solutions). But we can carry out the
Henneberg I construction and {\it choose} good edge lengths at
each step, which will guarantee that all subsequent steps yield
real solution. To do so, we will simply choose the edge lengths of
the two new edges to be approximately equal, and sum up to a bit
over the maximum distance (over all embeddings of the Henneberg I
graph $G$ to which this step is being applied) between the
endpoints in $G$ of the new edges. Hence:

\begin{prop}
Any realizable Henneberg I graph has at most $2^{n-2}$ embeddings.
For any Henneberg I graph, there exists a choice of edge lengths
such that the resulting framework has exactly $2^{n-2}$
embeddings.
\end{prop}

\medskip
\noindent {\bf Number of Embeddings for the Desargues Framework.}
To beat this bound we must use Henneberg II graphs. For special
edge lengths, the Desargues graph gives the smallest example with
more than $2^{n-2}$ embeddings. Namely, on $6$ vertices, it can
have up to $24$ ($> 2^{6-2}=16$) embeddings.

\begin{lemma}
\label{des} There exist edge lengths for the Desargues graph which
induce $24$ embeddings.
\end{lemma}

First, notice that if you remove an edge from a Laman graph, the
resulting graph is a one-degree-of-freedom mechanism. Once an edge
is pinned down, the other vertices trace algebraic curves. The
curve of the largest degree could then be intersected with a
circle whose center is placed at another vertex to get
intersection points allowing to place a new bar between the center
of the circle and one of these intersection points. Notice also
that such  mechanisms contain several rigid subcomponents: the
distances between vertices inside one such component stay fixed,
no matter how the mechanism is moving.

Henneberg II graphs must contain at least one vertex of degree
$3$. If that vertex (and its adjacent edges) are removed, the
resulting graph is a mechanism, whose vertices trace curves. The
intersection of these curves with circles should (in general) be
relatively easy to compute (even graphically, using for instance
the Cinderella \cite{cinderella} software). To be able to easily
put the three edges back in a way that would give as many
embeddings as possible for the given edge lengths, it would help
if two of these edges would stay at fixed distance (i.e. have an
edge between them, or belong to a rigid component of the
mechanism, which we will assume is grounded when the curve is
traced). To achieve this goal, we use a {\it floating center}
idea: place a circle anywhere, of any radius, as long as it gives
the desired number of crossings with the curve traced by the
mechanism. Then join the center via two bars to the grounded bar
of the mechanism: this yields the grounded triangle of the
Desargues configuration. To complete the construction, add a bar
between the center of the circle and any of the intersection
points. This can be done in as many ways as we had crossings. An
additional set of embeddings is obtained by flipping (about the
grounded edge) the two edges that were used to ground the center
of the floating circle.

It turns out that the Desargues graph has all the properties to
make this construction work.

\begin{figure}[htbp]
\vspace{-0.3in}
\begin{center}
\ \psfig{file=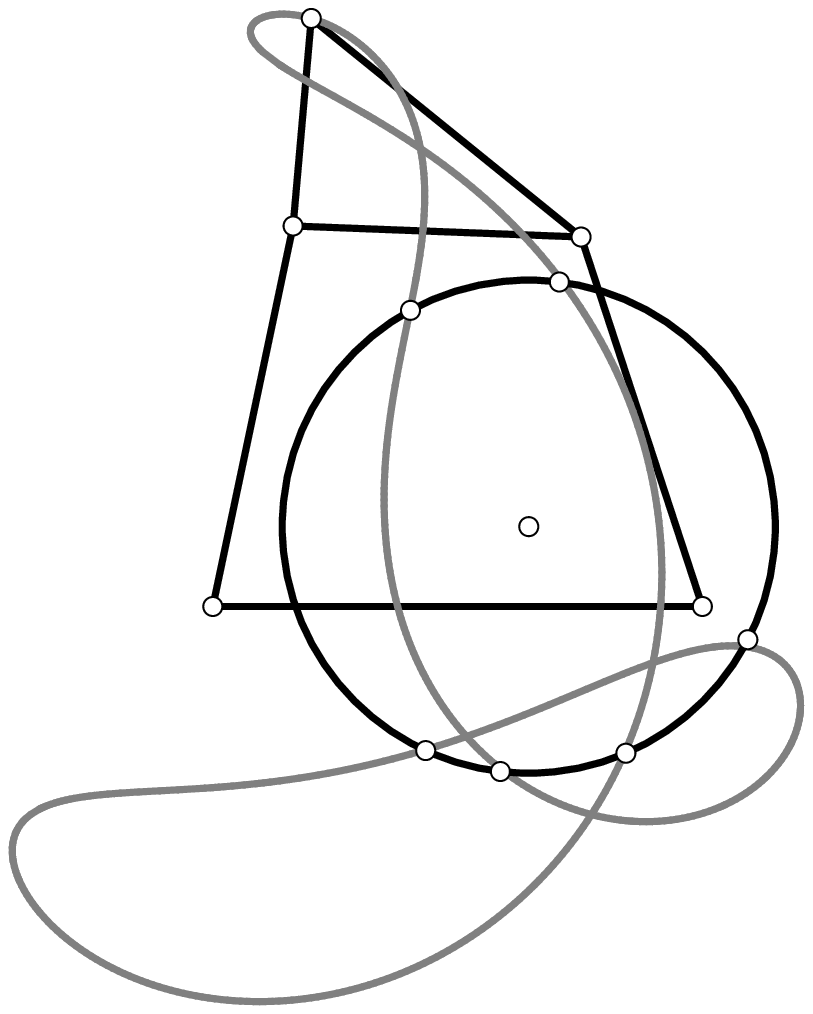,height=2.2in,width=1.8in}
\psfig{file=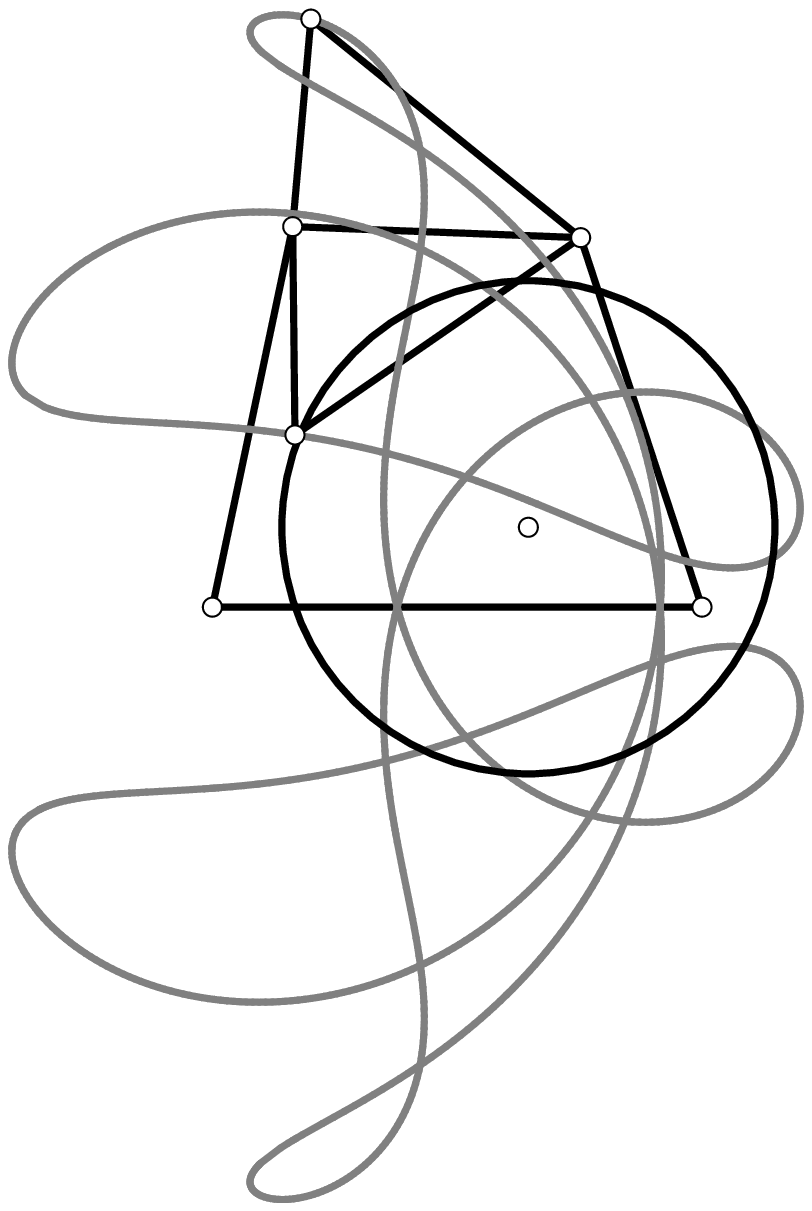,height=2.2in,width=1.8in}
\psfig{file=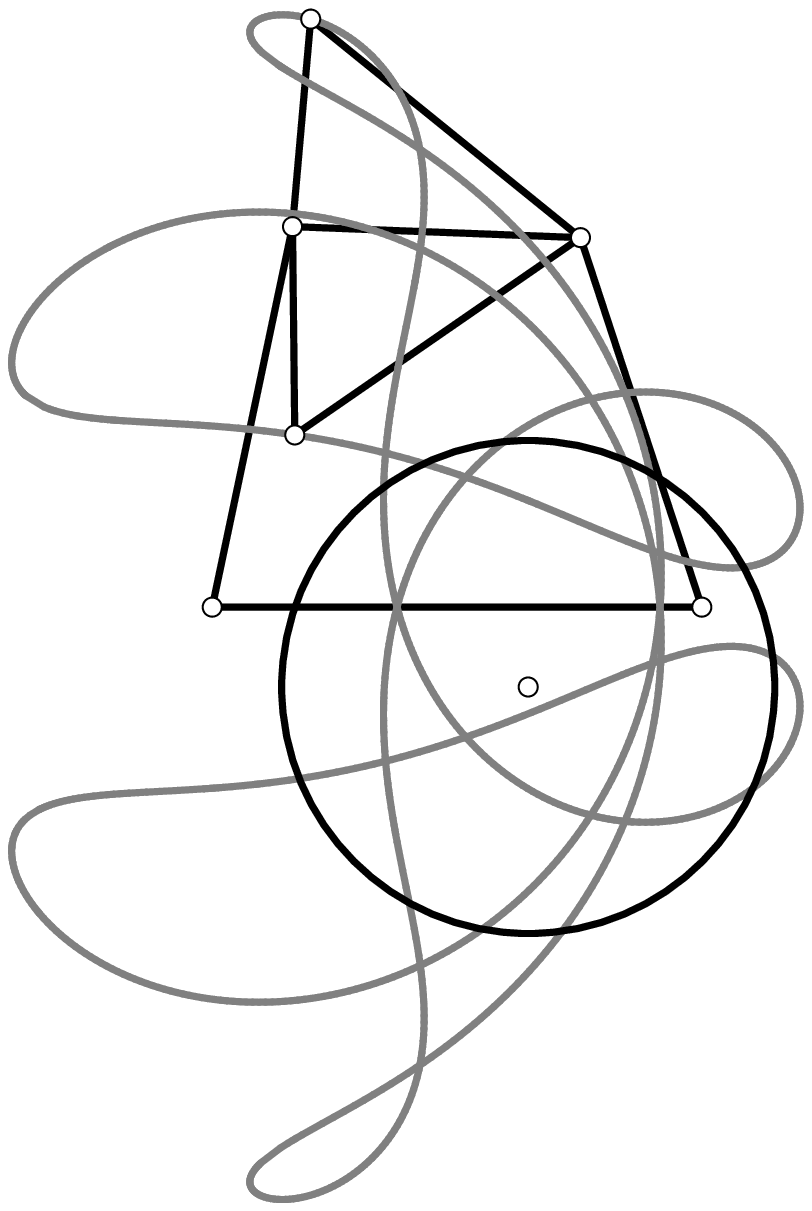,height=2.2in,width=1.8in}
\end{center}
\vspace{-0.3in}
\caption{
Left: one of the coupler curves, the ``floating'' circle and their
$6$ crossings. The bottom bar is grounded to the plane, the upper
three vertices are mobile and the top vertex traces the curve.
Middle: The $2$ coupler curves and a placement of a circle
crossing each of them in $6$ points. Right: The symmetric
placement of the circle from the middle picture induces another
$12$ crossings. } \figlab{des3a}
\end{figure}



The construction is illustrated in Figures \figref{des1} and
\figref{des3a}. Fig. \figref{des1} shows, in two embeddings,  the
mechanism obtained by removing a vertex of degree $3$. It is
assumed that the bottom edge is grounded and the motion of the
mechanism is guided by the rotation of one of the adjacent edges
around a grounded vertex. The curves traced by the degree $2$
vertex of the triangle attached to this $4$-bar mechanism have
received a lot of attention in mechanical engineering
(kinematics): indeed, in the '50's a whole atlas of such ``coupler
curves'' has been published (\cite{hrones}, see also
\cite{norton}). It is known that they are curves of degree $6$.

We have used Cinderella (\cite{cinderella}) to get a favorable
arrangement for the coupler curves. Fig. \figref{des3a}(left)
shows one curve and a position of a ``floating'' circle
intersecting it in $6$ points. The static black-and-white pictures
shown here do not capture the kinematic intricacy of the
arrangement. A Cinderella applet allowing the user to vary the
edge lengths and experiment with various shapes of the curves is
available from:{\texttt
http://cs.smith.edu/~streinu/Research/Embed/coupler.html}.

The two flipped triangles on the coupler edge give two distinct
and symmetric curves of degree $6$ each. We can find a position of
a circle intersecting each of them in $6$ points, see Figure
\figref{des3a}(middle). The symmetric arrangement gives another
$12$ crossings, for a total of $24$. See Fig.
\figref{des3a}(right). This proves Lemma \ref{des}.

\begin{figure}[htbp]
\vspace{-0.2in}
\begin{center}
\ \psfig{file=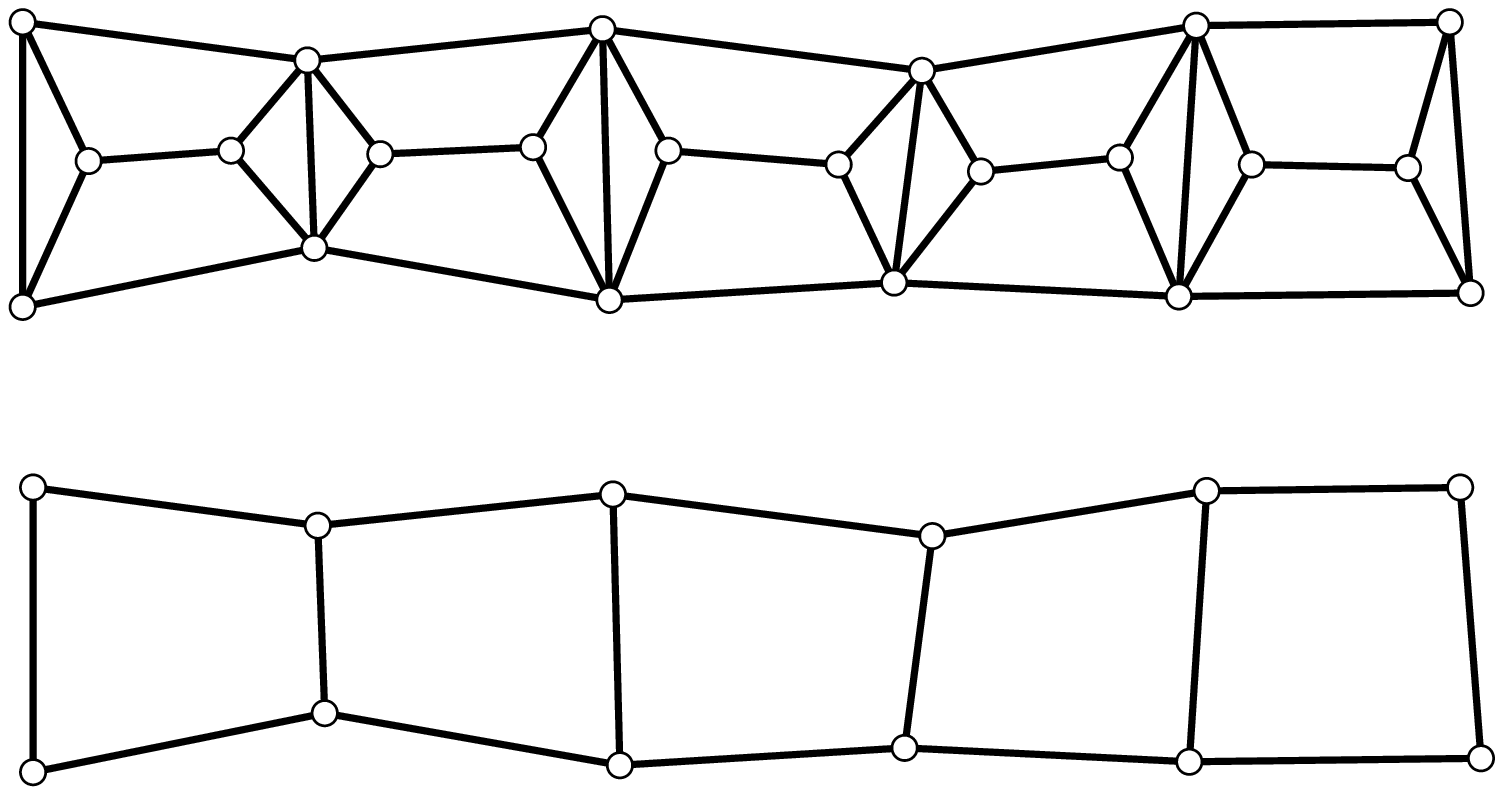,height=2.0in}
\psfig{file=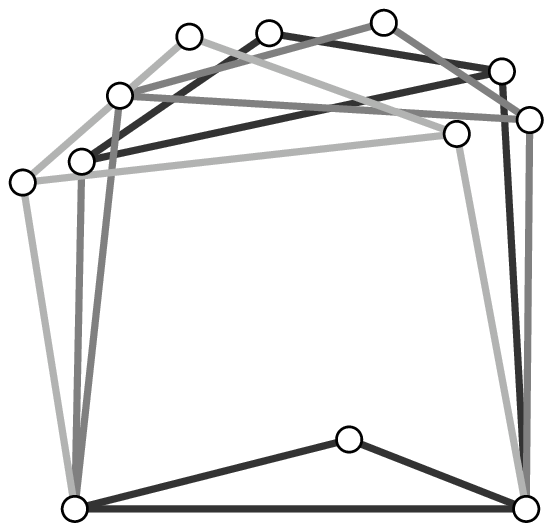,height=2.3in}
\end{center}
\vspace{-0.2in}
\caption{
Left: The ``caterpillar'' construction and its underlying
structure. Right: The structure of the ``fan'' construction.}
\figlab{catterpb}
\end{figure}



\medskip
\noindent {\bf Iterating the Desargues framework.} A simple
observation for an iterative construction is that if we glue
together, along any edge, two Desargues frameworks, we get a Laman
graph. For instance, we could glue them in a ``caterpillar''
fashion, as in Fig. \figref{catterpb}(left). The structure of this
construction is illustrated in Fig. \figref{catterpb}(middle): it
is a planar graph with quadrilateral faces (obtained from the
Desargues graphs by retaining just one of the underlying $4$-bar
mechanisms used in obtaining many embeddings) and whose dual graph
is a path. Any such glueing of quadrilaterals whose dual is a tree
would give the same bound. Intuitively, this is the equivalent of
glueing together triangles to get a polygon triangulation (the
dual of a polygon triangulation being a tree).

\begin{lemma}
\label{bound1} There exist edge length for which the number of
embeddings of the iterated ``caterpillar'' Desargues frameworks is
of the order of $24^{\frac n 4} \approx (2.21)^n$.
\end{lemma}

\begin{proof}
Each Desargues graph adds $4$ new vertices and multiplies by $24$
the total number of possible embeddings.
\end{proof}

To get a better bound we use the following observation.

\begin{lemma}
\label{perturb} A small perturbation of the lengths of the moving
edges of the $4$-bar mechanism from Lemma \ref{des} yields
slightly perturbed coupler curves, which still intersect the two
symmetric circles in $24$ points.
\end{lemma}

We now fix a bar in the plane, place on it a $4$-bar mechanism
satisfying Lemma \ref{des}, fix the positions of the symmetric
triangles and the corresponding circles that give the $24$
embeddings. Then perturb the $4$-bar mechanism (but not the base
triangle) several times as in Lemma \ref{perturb}, to obtain a
{\it fan}-like glueing of Desargues configurations. See Figure
\figref{catterpb}(right). Then the same two symmetric circles (or
circles of slightly perturbed radii, but with the same center)
cross each pair of coupler curves (one pair for each of these
mechanisms) in $24$ points, giving the desired number of
embeddings.

\begin{prop}
\label{bound2} There exist edge lengths for which the number of
embeddings of the iterated ``fan'' Desargues frameworks is of the
order of $24^{\frac n 3} \approx (2.88)^n$.
\end{prop}

\begin{proof}
Each Desargues graph adds now $3$ new vertices and $24$ possible
embeddings.
\end{proof}

As a final comment, we note that we can combine the perturbation
idea from the previous construction with the ``caterpillar
structure'' glueing to get another family of examples achieving
the same bound.

\medskip
\noindent {\bf Acknowledgements.} We thank Brigitte Servatius and
Walter Whiteley for pointing out useful references.


\bibliographystyle{abbrv}

\end{document}